\documentclass{amsart}
\usepackage{amssymb}
\usepackage{latexsym}
\usepackage{amscd}
\newtheorem{thm}{Theorem}
\newtheorem{tom}{Theorem}
\newtheorem{lem}{Lemma}[subsection]
\newtheorem{cor}{Corollary}
\newtheorem{coro}{Corollary}[subsection]
\newtheorem{crl}{Corollary}

\newtheorem{prop}{Proposition}
\newtheorem{propo}{Proposition}[subsection]
\newtheorem{prn}{Proposition}
\newtheorem{pd}{Proposition-Definition}[subsection]

\begin{document}
\renewcommand{\abstractname}{Abstract}
\begin{abstract}
Consider a $2$-plane $P \subset \mathbb{C}^n$ and let $D$ be a bounded region in $P$ with a piecewise-smooth boundary. Let $I(D)$ be the infimum of areas of all piecewise-smooth isotropic surfaces in $\mathbb{C}^n$ with the same boundary as $D$. Then $I(D)= \lambda_P^n \cdot Area(D)$. If $P$ is not complex, $\lambda_P^n < \frac{3\pi}{2\sqrt{2}}$. For a complex plane $\mathbb{C}  \subset \mathbb{C}^n$, $\lambda_{\mathbb{C}}^n \geq 2$, $\lambda_{\mathbb{C}}^2 \geq 3$ and also 
$\frac{3\pi^2}{2\sqrt{2}}$  is the area of an explicit Hamiltonian stationary isotropic Mobius band embedded in $\mathbb{C}^n$ whose boundary is a unit circle in $\mathbb{C}$. \\
As a corollary, a compact surface $\Sigma$ (possibly with boundary) in a symplectic manifold can be approximated by isotropic surfaces of area $\leq  \frac{3\pi}{2\sqrt{2}} Area(\Sigma)$. Another corollary is that a closed curve of length $l$ in $\mathbb{C}^n$ bounds an isotropic surface of area $\leq \frac{3l^2}{8\sqrt{2}}$. A related result is the following: consider $\mathbb{C}P^1 \subset \mathbb{C}P^n$ and let $D$ be a region in $\mathbb{C}P^1$. Let $I(D)$ be the infimum 
of areas of all isotropic surfaces in $\mathbb{C}P^n$ with the same boundary as $D$ representing the same relative homology class mod $2$ as $D$. Then $ 2 \cdot Area(D) \leq I(D) \leq \lambda_{\mathbb{C}}^n \cdot Area(D)$. Moreover the first inequality becomes an equality for $D=\mathbb{C}P^1$. 
\end{abstract}   
\title[Area comparison...]{Area comparison results for isotropic surfaces}     
\authors{Edward Goldstein\footnote{This material is based upon work supported by the National Science Foundation under agreement No. DMS-0111298. Any opinions, findings and conclusions or recommendations expressed in this material are those of the author and do not necessarily reflect the views of the National Science Foundation.}}
\maketitle
\section{Introduction}
This paper is concerned with area comparison results for isotropic surfaces in symplectic manifolds - those are surfaces on which the symplectic form restricts to $0$. The initial motivation comes from a paper \cite{Qiu} of  W. Qiu where it was shown that given a closed curve $\gamma \subset \mathbb{C}^n$ it bounds an isotropic surface $S$ with $Area(S) \leq C \cdot length(\gamma)^2$ for a constant $C$. This surface is necessarily non-orientable if the integral of the primitive $xdy$ of 
the symplectic form $dx \wedge dy$ over $\gamma$ is non-zero. One corollary of this is that isotropic surfaces are dense among all surfaces in the flat norm topology. It would be 
interesting to understand what is the infimum of areas of all isotropic 
surfaces with a given boundary $\gamma$. One can pose an analogous question 
for a null-homologous curve $\gamma$ in a K\"ahler manifold - in this case 
one should restrict attention to isotropic surfaces in some relative 
homology class. A related question is how well in terms of area comparison can one approximate an arbitrary surface by isotropic ones.\\
The first step in addressing those issues would be for isotropic surfaces bounding planar curves in $\mathbb{C}^n$. In this paper we establish the following theorem:
\begin{tom}
\label{maintheorem}
Consider a $2$-plane $P \subset \mathbb{C}^n$ and let $D$ be a bounded 
region in $P$ with a piecewise smooth boundary. Let $I(D)$ be the 
infimum of areas of all isotropic surfaces in $\mathbb{C}^n$ with the same boundary as $D$. Then $I(D)= \lambda_P^n \cdot Area(D)$.\\
 If $P$ is not complex, $\lambda_P^n < \frac{3\pi}{2\sqrt{2}}$. For a complex plane $\mathbb{C}  \subset \mathbb{C}^n$, $\lambda_{\mathbb{C}}^n \geq 2$, $\lambda_{\mathbb{C}}^2 \geq 3$ and also 
$\frac{3\pi^2}{2\sqrt{2}}$  is the area of an explicit Hamiltonian stationary isotropic Mobius band embedded in $\mathbb{C}^n$ whose boundary is a unit circle in $\mathbb{C}$.
\end{tom}
This yields the following corollaries:
\begin{cor}
\label{cor1}
Given a compact surface $\Sigma$ (possibly with boundary) in a symplectic manifold $(M,\omega,J)$, there is a sequence $S_n$ of isotropic surfaces in $M$ with $\partial(S_n)=\partial(\Sigma)$, $S_n \rightarrow \Sigma$ and $~ limsup ~Area(S_n) \leq  \frac{3\pi}{2\sqrt{2}} Area(\Sigma)$.
\end{cor}
Here $S_i$ converge to $\Sigma$ both in the flat norm topology and in the distance topologies. Combining the above corollary with the classical isoperimetric inequality in the Euclidean space gives
\begin{cor} 
\label{cor2}
Let $C$ be a closed curve in $\mathbb{C}^n$ of length $l$. Then $C$ bounds an isotropic surface $S \subset \mathbb{C}^n$ with $Area(S) \leq \frac{3l^2}{8\sqrt{2}}$.
\end{cor}
{\bf Remark:} The optimality of the estimates in the previous two corollaries hinges upon the question whether $\lambda_{\mathbb{C}}^n$ is indeed $\frac{3\pi}{2\sqrt{2}}$. The lower bound on $\lambda_{\mathbb{C}}^2$ says that in any case the estimates are pretty sharp, at least in dimension $4$.\\
A related result is the following
\begin{prop}
\label{projective}
Consider $\mathbb{C}P^1 \subset \mathbb{C}P^n$ and let $D$ be a region in 
$\mathbb{C}P^1$ with a piecewise smooth boundary. Let $I(D)$ be the infimum 
of areas of all isotropic surfaces in $\mathbb{C}P^n$ with the same boundary as $D$ representing the same relative homology class mod $2$ as $D$. Then $ 2 \cdot Area(D) \leq I(D) \leq \lambda_{\mathbb{C}}^n \cdot Area(D)$. Moreover the first inequality becomes an equality for $D=\mathbb{C}P^1$. 
\end{prop}
{\bf Remark:} For various regions $D \subset \mathbb{C}P^1$, $I(D)/Area(D)$ interpolates between $2$ ($D=\mathbb{C}P^1$) and $\lambda_{\mathbb{C}}^n$ (for $D$ concentrating near a point). It would be interesting to determine $I(D)$ for $D$ being a hemisphere in $\mathbb{C}P^1$.\\ 
\\
The paper is organized as follows: In section \ref{lambdapn} we'll show that $\lambda_P^n$ are well-defined. In section \ref{mobius} we'll use various isotropic Mobius bands to give upper bounds on $\lambda_P^n$. In section \ref{int} we'll present an integral-geometric formula due to R. Howard \cite{How} and also derive a new formula for areas of surfaces in $\mathbb{C}^n$, extending Howard's results. In section \ref{cp} we'll prove the lower bound in Proposition \ref{projective} and also find the ratio between the complex and the isotropic angles. In section \ref{c2} we'll give the lower bounds on $\lambda_{\mathbb{C}}^n$ and $\lambda_{\mathbb{C}}^2$, completing the proof of our main theorem. In section \ref{corol} we'll establish Corollaries \ref{cor1} and \ref{cor2} and Proposition \ref{projective}.\\
\\
{\bf Acknowledgements:} The author would like to express his gratitude to Richard Schoen for numerous discourses on Lagrangian surfaces. Special thanks go to Aleks, Olga and Nelly Neimark for their hospitality during the author's sojourn in Princeton, NJ.
\section{Definition of $\lambda_P^n$}
\label{lambdapn}
In this section we'll define $\lambda_P^n$ for a $2$-plane $P \subset \mathbb{C}^n$. First recall the result W. Qiu:
\begin{propo}
\label{qiu1}
\cite{Qiu} Given a closed curve $\gamma$ is $\mathbb{C}^n$ it bounds an isotropic 
surface $S$ whose area is $\leq C \cdot length(\gamma)^2$.
\end{propo}
This proposition has an immediate corollary:
\begin{coro}
\label{approx}
Given a bounded region $D$ on the $2$-plane  $P \subset 
\mathbb{C}^n$ with a piecewise smooth boundary there is an isotropic surface 
$S$ with the same boundary as $D$ whose area is $\leq 17C Area(D)$
\end{coro}
{\bf Proof:} Cover $P$ by a mesh of squares of size 
$\rightarrow 0$. This grid splits $D$ into regions $D_i$, most of them 
squares. Apply Qiu's result to each such region. Q.E.D.\\
This enables us to define $\lambda_P^n$:
\begin{lem}
\label{approx2}
Let $D_1$ be a unit square in $P$ and let $\lambda_P^n=I(D_1)$. For any bounded region $D$ in $P$ with a piecewise smooth boundary, \[I(D) =\lambda_P^n \cdot Area(D) \] 
\end{lem}
{\bf Proof:} Given an arbitrary $D$  we split the plane 
into a mesh of squares of size $\epsilon$. This splits $D$ into regions 
$D_i$. We find isotropic surfaces $S_i$ with the same boundary as $D_i$ as follows:\\
a) For those $D_i$ which are squares we require that $Area(S_i)$ is close to $\lambda_P^n \cdot Area (D_i)$.\\
b) For those $D_i$ which are not squares we use Corollary \ref{approx}.\\
We build a surface $S$ to be the union of $S_i$ to conclude that
\[ I(D) \leq \lambda_P^n \cdot Area(D) \]
To show that $I(D) = \lambda \cdot Area(D)$ we scale $D$ to be contained in a unit square $D_1$. We have that \[\lambda_P^n= I(D_1) \leq I(D)+I(D_1 - D) \leq \lambda_P^n(Area(D)+Area(D_1-D))\]
and all the inequalities become an equality. Q.E.D.
\section{Isotropic Mobius bands and upper bounds on $\lambda_P^n$}
\label{mobius}
\subsection{Construction of the Mobius bands.}
\label{mobi}
In this section we'll show that $\lambda_P^n \leq \frac{3\pi}{2\sqrt{2}}$ with a strict inequality for $P$ non-complex. We first exhibit explicit examples of isotropic Mobius bands with boundary in $P$, slightly generalizing the construction of D. Allcock \cite{All} and W. Qiu \cite{Qiu}. 
\begin{propo}
\label{dan}
Let $P$ and $P'$ be two $2$-planes in $\mathbb{C}^n$ which are $\omega$-orthogonal (here $\omega$ is the standard symplectic form on $\mathbb{C}^n$).
Let $\alpha(t), \beta(t)$ be two curves in $P$ and $P'$ correspondingly with
$\omega(\alpha,\alpha')= \omega(\beta,\beta')$. Consider a surface \[F(t,s)=\cos s 
\alpha(t)+ \sin s \beta(t)\] in $\mathbb{C}^n$. Then it is isotropic.
\end{propo}
{\bf Remark:} For $P$ and $P'$ being complex, orthogonal planes this proposition appears in \cite{All}.\\
Now we can think of the Mobius band as a rectangle $[0,2\pi] \times 
[0,\pi/2]$ on the $t,s$-plane with the identifications $(0,s) \simeq 
(2\pi,s)$ and $(t,\pi/2) \simeq (t+\pi,\pi/2)$. If we choose periodic 
$\alpha(t),\beta(t): [0,2\pi] \rightarrow \mathbb{C}$ with $\beta(t)=\beta(t+\pi)$ then the map $F(t,s)$ 
as in Proposition \ref{dan} will give an isotropic Mobius band whose boundary is $\alpha$.
\subsection{The case of a complex plane}
\label{compmobi}
We think of $P$ as the plane $z_2=0$ in $\mathbb{C}^2 \subset \mathbb{C}^n$ and $P'$ is the plane $z_1=0$ in $\mathbb{C}^2 \subset \mathbb{C}^n$. We 
choose $\alpha(t)=(e^{it},0)$ and $\beta(t)=(0,\frac{1}{\sqrt{2}}e^{2it})$. We 
compute the partial derivatives \[F_t=(i\cos s \cdot e^{it}, \sqrt{2}i 
\sin s \cdot e^{i2t}) ~ , ~ F_s= (-\sin s \cdot e^{it}, \frac{\cos s}{\sqrt{2}} \cdot e^{2it}) \]  
We note that $F_s$ and $F_t$ are orthogonal and 
\begin{equation}
\label{length}
|F_t|=\sqrt{2} |F_s|=\sqrt{1+\sin^2 s} 
\end{equation}
\begin{equation}
\label{mob}
Area(F)=\int_{t=0}^{2\pi}\int_{s=0}^{\pi/2}\frac{1+\sin^2s}{\sqrt{2}}dsdt=\frac{3}{2\sqrt{2}}\pi^2
\end{equation}
\begin{propo}
The surface $F$ is critical for the Area functional among isotropic surfaces with fixed boundary. 
\end{propo}
{\bf Proof:} Let $L_t$ be a $1$-parameter family of such surfaces with $L_0=F$. The deformation $\frac{d}{dt}L_t$ is realized by a vector field $v$ along $F$ which vanishes on the boundary. Also the isotropic condition implies that the contraction $i_v \omega$ of the symplectic form $\omega$ by $v$ is a closed $1$-form on $F$. Since this form vanishes on the boundary of $F$ and this boundary is a generator of the first homology of $F$ we get that $i_v \omega|_F$ is an {\it exact} $1$-form on $F$ i.e. $i_v \omega|_F=df$ for a function $f$ on $F$. Moreover we can choose $f$ to vanish on the boundary of $F$. Let $H$ be the mean curvature vector of $F$ and let $\sigma=i_H \omega|_F$ be the mean curvature $1$-form of $F$. The first variation formula tells that 
\begin{equation}
\label{fv1}
\frac{d}{dt} Area (L_t) =-\int_{F} v \cdot H= -\int_{F}df \cdot \sigma
\end{equation}
We use partitions of unity to write $f =\Sigma f_i$ there each $f_i$ is supported on an orientable piece of $F$ and vanishes on the boundary. Thus
\begin{equation}
\label{fv2}
\frac{d}{dt} Area (L_t) =-\Sigma \int_{F}df_i \wedge \ast \sigma
\end{equation}
Here $\ast \sigma$ is the Hodge star of the mean curvature $1$-form $\sigma$ (in the corresponding orientation on the support of $f_i$). One computes $\sigma$ as follows: we have a holomorphic $(2,0)$-form $dz_1 
\wedge dz_2$ on $\mathbb{C}^2$. When we restrict it to the Lagrangian 
surface $F$ we have that \[dz_1 \wedge dz_2|_F = e^{i\theta} vol(F)\] 
Here $vol(F)$ is the area form of $F$, defined up to a sign, and the mean 
curvature 1-form  is $\sigma= d\theta$ (see \cite{Oh}, Proposition 2.2 or \cite{CG}). Using the tangent 
vectors $F_t$ and $F_s$ we see that $\sigma=3dt$ and its Hodge star in the corresponding orientation is $\ast \sigma= \frac{3}{\sqrt{2}} ds$ and it is also closed. Such a surface is called Hamiltonian stationary - see \cite{ScW}. We get that the first variation of area in equation (\ref{fv2}) vanishes.   Q.E.D. 
\subsection{An upper bound on $\lambda_P^n$.}
\label{upperbound}
We now study the case of a non-complex plane $P \subset \mathbb{C}^n$. We can assume w.l.o.g. that $P=~ span ~((1,0),(ia,\sqrt{1-a^2})) \subset \mathbb{C}^2 \subset \mathbb{C}^n$ for $0 \leq a < 1$. We can take $P'=~ span ~((0,1),(\sqrt{1-a^2},ia)) \subset \mathbb{C}^2 \subset \mathbb{C}^n$. We can apply Proposition \ref{dan} with the curves \[\alpha(t)=\cos t(1,0)+\sin t (ia,\sqrt{1-a^2}) ~ , ~ \beta(t)= \frac{\cos2t}{\sqrt{2}}(0,1)+\frac{\sin 2t}{\sqrt{2}}(\sqrt{1-a^2},ia)\]
We get an isotropic Mobius band as in section \ref{mobi}. We compute the tangent vectors
\[F_t=\cos s(-\sin t(1,0)+\cos t (ia,\sqrt{1-a^2})) +\sqrt{2}\sin s(-\sin 2t(0,1)+ \cos 2t(\sqrt{1-a^2},ia))\]
\[F_s=-\sin s(\cos t(1,0)+\sin t (ia,\sqrt{1-a^2})) +\frac{\cos s}{\sqrt{2}}(\cos 2t(0,1)+ \sin 2t(\sqrt{1-a^2},ia))\]
\[F_t \cdot F_t=\cos^2 s +2\sin^2s-2\sqrt{2(1-a^2)} \cos s \sin s \sin 3t ~ , ~ F_s \cdot F_s =\frac{F_t \cdot F_t}{2}\]
We note that the area of the parallelogram spanned by two vectors $X$ and $Y$ is $\sqrt{(X\cdot X )(Y \cdot Y) -(X \cdot Y)^2} \leq \sqrt{(X\cdot X )(Y \cdot Y)}$. We conclude that
\[Area(F) < \int_{s=0}^{\pi/2} \int_{t=0}^{2\pi} \frac{F_t \cdot F_t}{\sqrt{2}} d s d t=  \int_{s=0}^{\pi/2} \int_{t=0}^{2\pi}\frac{1 +\sin^2s}{\sqrt{2}}- \sin 2s \sin 3t \sqrt{1-a^2} ds dt =\frac{3\pi^2}{2\sqrt{2}}\]
Thus for $P$ non-complex, we got an isotropic Mobius band $F$ whose boundary is a unit circle in $P$ with $Area(F) < \frac{3\pi^2}{2\sqrt{2}}$. Thus $\lambda_P^n <  \frac{3\pi}{2\sqrt{2}}$.
\section{Formulas from integral geometry}
\label{int}
\subsection{Intersections of compact submanifolds}
In this section we'll describe a formula from integral geometry following the exposition in R. Howard \cite{How}. Let $G$ be a unimodular Lie group (i.e. it admits a bi-invariant volume form) and let $K \subset G$ be a compact subgroup. Pick a left invariant metric $g$ on $G$ which is also right invariant under $K$. Let $M=G/K$ be the corresponding homogeneous space. Let $P$ and $Q$ be submanifolds of $M$ of complementary dimensions. For a point $p \in P$ and $q \in Q$ we define an angle $\sigma(p,q)$ between the tangent planes $T_pP$ and $T_qQ$ as follows: First we choose some elements $g$ and $h$ in $G$ which move $p$ and $q$ respectively to the same point $r \in M$. Now the tangent planes $g_{\ast}T_pP$ and $h_{\ast}T_qQ$ 
are in the same tangent space $T_r M$ and we can define an 
angle between them as follows: take an orthonormal basis $u_1,\ldots,u_k$ for $g_{\ast}T_pP$ and an orthonormal basis $v_1 \ldots v_l$ for 
$h_{\ast}T_qQ$ and define \[\sigma(g_{\ast}T_pP,h_{\ast}T_qQ)= |u_1 \wedge 
\ldots \wedge v_l|\]
The later quantity $\sigma(g_{\ast}T_pP,h_{\ast}T_qQ)$ depends on the 
choices $g$ and $h$ we made. To mend this we'll need to average this out by the stabilizer group $K$ of the point $r$. Thus we define:
\begin{equation}
\label{angle}
\sigma(p,q)=\int_{K} \sigma(g_{\ast}T_pP, k_{\ast}h_{\ast}T_qQ) dk
\end{equation}
Now assume that $P$ and $Q$ are compact, possibly with boundary. There is a following general formula due to R. Howard \cite{How} :
\begin{equation}
\label{main}
 \int_{G} \#(P \bigcap gQ) dg= \int_{P \times Q} \sigma(p,q) dp dq
\end{equation}
\subsection{Intersections with complex hyperplanes}
In this section we'll derive a new formula for areas of surfaces in $\mathbb{C}^n$, which is a corollary of equation \ref{main}. Let $G$ be the group of isomorphisms (biholomorphic isometries) of $\mathbb{C}^n$. Then $G$ is the semidirect product of $U(n)$ with $\mathbb{C}^n$. Here $\mathbb{C}^n$ acts on itself by translations and $K=U(n)$ is the stabilizer of the origin in $\mathbb{C}^n$. Also $[U(n),\mathbb{C}^n] \subset \mathbb{C}^n$. Hence $G$ is unimodular.\\
Let $N$ be the space of all complex hyperplanes in $\mathbb{C}^n$ (not necessarily passing through the origin). Then $N$ is a homogeneous space $N=G/H$ there $H$ is a stabilizer of a hyperplane (say $\mathbb{C}^{n-1}$) in  $\mathbb{C}^n$. On shows as before that $H$ is also unimodular and the space $N=G/H$ has a $G$-invariant volume form. \\
Let $B$ be a unit ball in a hyperplane $\mathbb{C}^{n-1} \in N$ and the stabilizer of this hyperplane in $G$ is $H$ as before. Pick a point $q$ on $\mathbb{C}^{n-1}$ and define the number $A$ to be:
\begin{equation}
\label{unit}
A=vol (h \in H| q \in h(B))
\end{equation}
Clearly this number $A>0$ is independent of a choice of a point $ q \in \mathbb{C}^{n-1}$. We have the following lemma:
\begin{lem}
\label{new}
Let $P$ be a compact surface in  $\mathbb{C}^n$, possibly with boundary. 
For a point $p \in P$ we define $\sigma_{p,C}$ to be the angle which $T_pP$ forms with a complex hyperplane as in equation (\ref{angle}). Let $N$ be the space of all complex hyperplanes in $\mathbb{C}^n$. Then
\[\frac{A}{E_{n-1}}\int_{N} \#(\eta \bigcap P) d\eta= \int_{P} \sigma_{p,C} dp\]
Here $A$ is given by equation (\ref{unit}) and $E_{n-1}$ is the volume of a unit ball in $\mathbb{C}^{n-1}$.
\end{lem}
{\bf Proof:} Let $B$ be the unit ball in the complex hyperplane $\mathbb{C}^{n-1}$. We'll use equation (\ref{main}) with $Q=B$. We get that
\[ \int_{G} \#(P \bigcap gB) dg= E_{n-1} \int_{P} \sigma_{p,C} dp\]
We rewrite
\[ \int_{G} \#(P \bigcap gB) dg=\int_{N} (\int_{H} \#(P \bigcap \eta \cdot h (B))dh)d\eta\]
For a generic hyperplane $\eta \in N$, it intersects $P$ transversally in points $p_1,\ldots,p_l$. We have \[\int_{H} \#(P \bigcap \eta \cdot h (B))dh=\Sigma_{i=1}^{l}\int_{H} \#( \eta^{-1}(p_i) \bigcap h(B))dh=A \cdot l\]
and this proves the lemma. Q.E.D.
\section{Area comparison in $\mathbb{C}P^n$}
\label{cp}
In this case the group $G=SU(n+1)$ acts on $\mathbb{C}P^n$ with a stabilizer $K \simeq U(n)$. Thus we view $\mathbb{C}P^{n}=SU(n+1)/K$ and the Fubini-Study metric is induced from the bi-invariant metric on $SU(n+1)$. Let $P$ be a surface in $\mathbb{C}P^n$ and let $Q$ be a linear 
$\mathbb{C}P^{n-1} \subset \mathbb{C}P^n$. We'll treat two cases: $P$ is isotropic or complex.\\
{\bf Isotropic case:} 
Since $SU(n+1)$ acts transitively on the Grassmanian of isotropic planes 
in $\mathbb{C}P^n$ we conclude that this angle is a constant depending 
just on $n$:
\begin{equation}
\label{isotangle}
\sigma(p,q)=C_{I,n}
\end{equation}
{\bf Complex case:} Since $SU(n+1)$ acts transitively on the Grassmanian 
of complex planes in $\mathbb{C}P^n$ we conclude that this angle is a 
constant depending just on $n$:
\begin{equation}
\label{comangle}
\sigma(p,q)=C_{C,n}
\end{equation}
We'll use equations (\ref{isotangle},\ref{comangle},\ref{main}) for $P$ being the totally geodesic $\mathbb{R}P^2$ and $\mathbb{C}P^1$ correspondingly. In both cases for generic $g \in SU(n+1)$, $\#(P \bigcap gQ)=1$. Also $Area(\mathbb{R}P^2)=2\pi$, $Area(\mathbb{C}P^1)=\pi$. Hence we conclude that 
\begin{equation}
\label{ratio}
C_{C,n}=2C_{I,n}
\end{equation}
Now we can prove the lower bound in Proposition \ref{projective} stated in the introduction.
\begin{lem}
\label{CPN}
Consider $\mathbb{C}P^1 \subset \mathbb{C}P^n$ and let $D$ be a region in 
$\mathbb{C}P^1$ with piecewise smooth boundary. Let $I(D)$ be the infimum 
of areas of all isotropic surfaces in $\mathbb{C}P^n$ with the same boundary as $D$ representing the same relative homology class mod $2$ as $D$. Then $I(D) 
\geq 2 \cdot Area(D)$ with equality for $D=\mathbb{C}P^1$.
\end{lem}
{\bf Proof:} Let $S$ be an isotropic surface in the same relative homology 
class mod $2$ as $D$. We'll use formulas (\ref{ratio},\ref{main}) for $S$ and $D$. We note that for a generic $g \in SU(n+1)$, $\#(D 
\bigcap gQ)$ is either $0$ or $1$. If $\#(D \bigcap gQ)=1$ then the 
intersection number mod $2$ of $gQ$ with the relative homology class of 
$D$  is $1$. Hence $\#(S \bigcap gQ) \geq 1$. Since $C_{C,n}=2C_{I,n}$ we 
conclude that $Area(S) \geq 2 \cdot Area(D)$.\\
Also if $D=\mathbb{C}P^1$ then we can take $S=\mathbb{R}P^2$ which is in 
the same homology class mod $2$ as $D$ and $Area(S)=2 \cdot Area(D)$. 
Q.E.D.
\section{Lower bounds on $\lambda_{\mathbb{C}}^n$ and the main theorem}
\label{c2}
\subsection{Proof that $\lambda_{\mathbb{C}}^n \geq 2$}
The goal of this section is to give estimates for areas of isotropic 
surfaces in $\mathbb{C}^n$ whose boundary lies on a complex line 
$\mathbb{C}^1 \subset \mathbb{C}^n$. Let $G$ be the group of isomorphisms (biholomorphic isometries) of $\mathbb{C}^n$. Let $U(n) \subset G$ be the stabilizer of the origin in $\mathbb{C}^n$. Let $I_L$ be an angle between an isotropic plane and a complex hyperplane and let $I_C$ be an angle between a complex plane and a complex hyperplane as in equation (\ref{angle}).
\begin{propo}
\label{twice}
$I_C=2 \cdot I_L$
\end{propo}
{\bf Proof:} The result is equivalent to the corresponding equation (\ref{ratio}) in $\mathbb{C}P^n$. In both cases the stabilizer $K$ of a point is isomorphic to $U(n)$ with the standard action on the tangent space. We use the formula (\ref{angle}) to compute the angle. Since a left invatiant volume form on $U(n)$ is unique up to a constant multiple we conclude the statement of our proposition. Q.E.D.\\
\begin{propo}
\label{lambdan}
$\lambda_{\mathbb{C}}^n \geq 2$.
\end{propo}
{\bf Proof:} 
Let $D$ be a planar region in $\mathbb{C}^1\subset \mathbb{C}P^n$ and let $S$ be an isotropic surface in $\mathbb{C}^n$ with the same boundary as $D$. We'll use Lemma \ref{new} for $P= S$ and $P=D$. We have by Proposition \ref{twice}, $I_C=2 \cdot I_L$. Also any hyperplane $\eta$ that intersects $D$ transversally does so in exactly one point, hence the intersection number mod $2$ of $D$ and $\eta$ is $1$ and the same is true for $\eta$ and $S$, hence $\eta$ intersects $S$. Using Lemma \ref{new} we conclude that $Area(S) \geq 2 Area(D)$. Q.E.D.\\
The next sections show that in complex dimension $n=2$ one can say more.
\subsection{Refinement of Lemma \ref{new}}
Let $N=G/H$ be the space of all lines in $\mathbb{C}^2$ - here we use the notation of section \ref{int}. Let $\mathbb{C}P^1$ we the space of all lines in $\mathbb{C}^2$ passing through the origin. For any line $n \in N$ there is a unique line $\pi(n) \in \mathbb{C}P^1$ which is perpendicular to $n$. Thus we have a fibration \[\pi: N \mapsto \mathbb{C}P^1\] Now $G$ acts on $N$ and also on $\mathbb{C}P^1$ (one shifts the line back to the origin after acting on it by an element $g \in G$). This action commutes with  $\pi: N \mapsto \mathbb{C}P^1$ and it preserves the volume forms on $N$ and on $\mathbb{C}P^1$. Hence there are canonical volume forms $\mu$ on the fibers of $\pi$ such that $G$ acts by volume-preserving diffeomorphisms on the fibers.\\
For a line $\kappa \in \mathbb{C}P^1$, the fiber $\pi^{-1}(\kappa)$ is naturally identifined with the set of points on $\kappa$. Namely a point $x \in \kappa$ is identified with a line $\kappa(x)$ which passes through $x$ and perpendicular to $\kappa$. So
\begin{equation}
\label{fiber}
\pi^{-1}(\kappa)= \bigcup (\kappa(x)|x \in \kappa)
\end{equation}
Also the volume form $\mu$ on $\pi^{-1}(\kappa)$ equals to :
\begin{equation}
\label{lambda}
\mu= \delta vol(\kappa)
\end{equation}
Here $vol(\kappa)$ is the area form on $\kappa$ and $\delta$ is a constant independent of $\kappa$. If $P$ is a compact surface (possibly with boundary) then Lemma \ref{new} tells that
\begin{equation}
\label{integration}
\int_{P} \sigma_{p,C} dp=\frac{\delta A}{\pi} \int_{\mathbb{C}P^1} (\int_{\kappa} \#(\kappa(x) \bigcap P)dx )d \kappa
\end{equation}
To understand the later quantity $\int_{\kappa} \#(\kappa(x) \bigcap P) dx$ we make the following definitions: 
\begin{pd}
\label{def1}
For $\kappa \in \mathbb{C}P^1$ define a decomposable two-form $\rho_{\kappa}$ on $\mathbb{C}^2$ whose kernel is the line $\kappa^{\perp}$ perpendicular to $\kappa$ and such that $\rho_{\kappa}$ restricts to the  area form on $\kappa$. If $\kappa $ and $\kappa^{\perp}$ are perpendicular then \[\rho_{\kappa}+\rho_{\kappa^{\perp}}=\omega \]
Here $\omega$ is the K\"ahler form on $\mathbb{C}^2$.
\end{pd}
\begin{pd}
\label{replace}
For a compact surface $P \subset \mathbb{C}^2$ (possibly with boundary) and a line $\kappa$ passing through the origin in $\mathbb{C}^2$ define
\[ F(P,\kappa)=\int_{\kappa} \#(\kappa(x) \bigcap P) dx= \int_P |\rho_{\kappa}| dp\]
\end{pd}
Note that we can integrate the absolute value $|\rho_{\kappa}|$ even if $P$ is not orientable. 
Thus equation \ref{integration} means that for any compact surface $P$ (possibly with boundary) in $\mathbb{C}^2$,
\begin{equation}
\label{int2}
\int_{P} \sigma_{p,C} dp=\frac{\delta A}{\pi} \int_{\mathbb{C}P^1} F(P,\kappa) d \kappa
\end{equation}
We summarize the properties of $F(P,\kappa)$ in the following proposition:
\begin{propo}
\label{F}
Let $\kappa \in \mathbb{C}P^1$ be a line in $\mathbb{C}^2$ passing through the origin. If $S$ is a Lagrangian surface then $F(S,\kappa)=F(S,\kappa^{\perp})$. For a planar region $D \subset \mathbb{C}^1$, $F(D,\kappa)=Area(D) \cos^2 \alpha_{\kappa}$ there $\alpha_{\kappa}$ is an angle between $\kappa$ and $ \mathbb{C}^1$. Moreover if $S$ has the same boundary as $D$ then $F(S,\kappa) \geq F(D,\kappa)$.
\end{propo}
{\bf Proof:} We see from Proposition-Definition \ref{def1} that restricted to $S$, $|\rho_{\kappa}|=|\rho_{\kappa^{\perp}}|$. The first claim is now immediate from Proposition-Definition \ref{replace}. The second claim is obvious. For the third claim we note that the orthogonal projection of $S$ onto $\kappa$ contains the orthogonal projection of $D$ onto $\kappa$. Hence the claim follows from Proposition-Definition \ref{replace}. Q.E.D.
\subsection{Proof of the main theorem}
The last step in establishing Theorem \ref{maintheorem} is 
\begin{lem}
\label{low}
Consider $\mathbb{C}^1 \subset \mathbb{C}^2$ and let $D$ be a bounded 
region in $\mathbb{C}^1$ with a piecewise smooth boundary. Let $S$ be a piecewise smooth Lagrangian surface in $\mathbb{C}^2$ with the same boundary as $D$. Then 
$Area(S) \geq 3 \cdot Area(D) $.
\end{lem}
{\bf Proof:} We use equation (\ref{int2}) and Proposition \ref{F} for $P=D$ and $P=S$:
\begin{equation}
\label{AreaD}
Area(D)= Area(D) \cdot \frac{\delta A}{I_C \pi} \int_{\mathbb{C}P^1} \cos^2 \alpha_{\kappa} d \kappa
\end{equation}
\begin{equation}
\label{AreaP}
Area(S) \geq Area(D) \cdot \frac{\delta A}{I_L \pi} \cdot 2 \int_{\alpha_{\kappa} \leq \pi/4} \cos^2 \alpha_{\kappa} d \kappa
\end{equation}
So we need to compute the ratio between two numbers:
\[ Num_C=\int_{\mathbb{C}P^1} \cos^2 \alpha_{\kappa} d \kappa ~ , ~ Num_L=2 \int_{\alpha_{\kappa} \leq \pi/4} \cos^2 \alpha_{\kappa} d \kappa\]
Now $Num_C=1/2 \cdot Area(\mathbb{C}P^1)$. To understand $Num_L$ let $(z_1,z_2)$ be coordinates on $\mathbb{C}^2$ so that $\mathbb{C}^1=(z_2=0)$.
We introduce the inhomogeneous coordinate $z=x+iy=\frac{z_2}{z_1}$ on $\mathbb{C}P^1$. So \[\cos^2 \alpha_{\kappa}=\frac{1}{1+|z|^2} \]
The Fubini-Study form is
\[\omega_{FS}=i \partial \overline{\partial}ln(1+|z|^2)= \frac{2 dx \wedge dy}{(1+|z|^2)^2}\]
and the region $ \alpha_{\kappa} \leq \pi/4$ corresponds to $|z| \leq 1$.
We pass to polar coordinates and easily compute that 
\begin{equation}
\label{NLC} 
Num_L=1.5 \cdot Num_C 
\end{equation}
We combine equations (\ref{AreaD},\ref{AreaP},\ref{NLC}) and Proposition \ref{twice} to prove the lemma. Q.E.D.\\
We are now ready to establish the main theorem.
\begin{thm}
Consider a $2$-plane $P \subset \mathbb{C}^n$ and let $D$ be a bounded 
region in $P$ with a piecewise smooth boundary. Let $I(D)$ be the 
infimum of areas of all isotropic surfaces in $\mathbb{C}^n$ with the same boundary as $D$. Then $I(D)= \lambda_P^n \cdot Area(D)$.\\
If $P$ is not complex, $\lambda_P^n < \frac{3\pi}{2\sqrt{2}}$. For a complex plane $\mathbb{C}  \subset \mathbb{C}^n$, $\lambda_{\mathbb{C}}^n \geq 2$, $\lambda_{\mathbb{C}}^2 \geq 3$ and also 
$\frac{3\pi^2}{2\sqrt{2}}$  is the area of an explicit Hamiltonian stationary isotropic Mobius band embedded in $\mathbb{C}^n$ whose boundary is a unit circle in $\mathbb{C}$.
\end{thm}
{\bf Proof:} The fact that $\lambda_P^n$ is well defined was shown in section \ref{lambdapn}. Upper bounds on $\lambda_P^n$ were given in section \ref{upperbound}. Lower bound on $\lambda_{\mathbb{C}}^n$ is Proposition \ref{lambdan}. The lower bound on $\lambda_{\mathbb{C}}^2$ is Lemma \ref{low}. Finally the Hamiltonian stationary Mobius band was constructed in section \ref {compmobi}. Q.E.D. 
\section{Proof of the Corollaries and Proposition \ref{projective}}
\label{corol}
Now we can prove Corollaries \ref{cor1}, \ref{cor2} and Proposition \ref{projective}.
\begin{crl}
Given a compact surface $\Sigma$ (possibly with boundary) in a symplectic manifold $(M,\omega,J)$, there is a sequence $S_n$ of isotropic surfaces in $M$ with $\partial(S_n)=\partial(\Sigma)$, $S_n \rightarrow \Sigma$ and $~ limsup ~Area(S_n) \leq  \frac{3\pi}{2\sqrt{2}} Area(\Sigma)$.
\end{crl}
Here $S_i$ converge to $\Sigma$ both in the flat norm topology and in the distance topology. \\
{\bf Proof:} For any point $p \in \Sigma$ one has a symplectomorphism $\phi_p$ from a neighbourhood $U_p$ of $p$ in $M$ into $\mathbb{C}^n$. Moreover one can choose the differential $d\phi_p$ to be an isometry. Take $U_p$ small enough so that $\phi_p$ is $C^0$-close to an isometry. Covering $\Sigma$ by various $U_p$'s it is clearly enough to prove the statement for $\Sigma \subset \mathbb{C}^n$. But this follows from approximating $\Sigma$ by polygonal surfaces, the upper bound on $\lambda_P^n$ in Theorem \ref{maintheorem} and Proposition \ref{qiu1}. Q.E.D.
\begin{crl} 
Let $C$ be a closed curve in $\mathbb{C}^n$ of length $l$. Then $C$ bounds an isotropic surface $S \subset \mathbb{C}^n$ with $Area(S) \leq \frac{3l^2}{8\sqrt{2}}$.
\end{crl}
{\bf Proof:} The classical isoperimetric inequality says that $C$ bounds a surface $\Sigma$ in $\mathbb{C}^n$ of area $\leq \frac{l^2}{4\pi}$ with an equality iff $C$ is a circle on a plane. Apply Corollary \ref{cor1} to $\Sigma$. Q.E.D.
\begin{prn}
Consider $\mathbb{C}P^1 \subset \mathbb{C}P^n$ and let $D$ be a region in 
$\mathbb{C}P^1$ with a piecewise smooth boundary. Let $I(D)$ be the infimum 
of areas of all isotropic surfaces in $\mathbb{C}P^n$ with the same boundary as $D$ representing the same relative homology class mod $2$ as $D$. Then $ 2 \cdot Area(D) \leq I(D) \leq \lambda_{\mathbb{C}}^n \cdot Area(D)$. Moreover the first inequality becomes an equality for $D=\mathbb{C}P^1$. 
\end{prn}
{\bf Proof:} The lower bound and the statement for $D=\mathbb{C}P^1$ follows from Lemma \ref{CPN}. The upper bound follows from Corollary \ref{cor1}.  Q.E.D.

Brandeis University\\
egold@brandeis.edu

\end{document}